\documentclass[10pt]{article}

\usepackage{amsmath, amssymb, amscd, fancyhdr}

\newcommand{\version}{version 8.0,\ \ Oct. 04, 2010}

\def\eqref#1{(\ref{#1})}

\newcommand{\arrow}{{\:\longrightarrow\:}}
\newcommand{\Z}{{\Bbb Z}}
\newcommand{\C}{{\Bbb C}}
\newcommand{\R}{{\Bbb R}}

\newcommand{\6}{\partial}
\def\1{\sqrt{-1}\:}

\newcommand{\cntrct}                
{\hspace{2pt}\raisebox{1pt}{\text{$\lrcorner$}}\hspace{2pt}}

\def\Bbb#1{\mathbb #1}

\renewcommand{\bar}{\overline}
\renewcommand{\phi}{\varphi}
\renewcommand{\epsilon}{\varepsilon}

\newcommand{\End}{\operatorname{End}}

\newcommand{\Tw}{\operatorname{Tw}}

\renewcommand{\Re}{\operatorname{Re}}
\renewcommand{\Im}{\operatorname{Im}}

\newcommand{\comment}[1]{{}}

\def\blacksquare{\hbox{\vrule width 4pt height 4pt depth 0pt}}
\def\endproof{\blacksquare}

\makeatletter

\newcounter{Mycounter}[section]
\newcounter{lemma}[section]
\renewcommand{\thelemma}{\noindent{Lemma \thesection.\arabic{lemma}}}
\newcommand{\lemma}{%
     \setcounter{lemma}{\value{Mycounter}}
     \refstepcounter{lemma}
     \stepcounter{Mycounter}
     {\bf \thelemma:\ }}

\newcounter{claim}[section]
\renewcommand{\theclaim}{\noindent{Claim \thesection.\arabic{claim}}}
\newcommand{\claim}{%
     \setcounter{claim}{\value{Mycounter}}
     \refstepcounter{claim}
     \stepcounter{Mycounter}
     {\bf \theclaim:\ }}

\newcounter{sublemma}[section]
\newcounter{corollary}[section]
\renewcommand{\thecorollary}{\noindent{Corollary \thesection.\arabic{corollary}}}
\newcommand{\corollary}{%
     \setcounter{corollary}{\value{Mycounter}}
     \refstepcounter{corollary}
     \stepcounter{Mycounter}
     {\bf \thecorollary:\ }}

\newcounter{theorem}[section]
\renewcommand{\thetheorem}{\noindent{Theorem \thesection.\arabic{theorem}}}
\newcommand{\theorem}{%
     \setcounter{theorem}{\value{Mycounter}}
     \refstepcounter{theorem}
     \stepcounter{Mycounter}
     {\bf \thetheorem:\ }}

\newcounter{conjecture}[section]
\newcounter{proposition}[section]
\renewcommand{\theproposition}
       {\noindent{Proposition \thesection.\arabic{proposition}}}
\newcommand{\proposition}{%
     \setcounter{proposition}{\value{Mycounter}}
     \refstepcounter{proposition}
     \stepcounter{Mycounter}
     {\bf \theproposition:\ }}

\newcounter{definition}[section]
\renewcommand{\thedefinition}
       {\noindent{Definition~\thesection.\arabic{definition}}}
\newcommand{\definition}{%
     \setcounter{definition}{\value{Mycounter}}
     \refstepcounter{definition}
     \stepcounter{Mycounter}
     {\bf \thedefinition:\ }}

\newcounter{example}[section]
\newcounter{remark}[section]
\renewcommand{\theremark}{\noindent{Remark \thesection.\arabic{remark}}}
\newcommand{\remark}{%
     \setcounter{remark}{\value{Mycounter}}
     \refstepcounter{remark}
     \stepcounter{Mycounter}
     {\bf \theremark:\ }}

\newcounter{problem}[section]
\newcounter{question}[section]
\@addtoreset{equation}{section}
\@addtoreset{footnote}{section}
\makeatother

\begin{document}

\begin{center}
{\LARGE\bf Hodge theory on
nearly K\"ahler manifolds}
\\[4mm]
Misha Verbitsky\footnote{Misha Verbitsky is 
Partially supported by RFBR grant
10-01-93113-NCNIL-a, 
AG Laboratory SU-HSE, RF government grant, ag. 11.G34.31.0023, 
Science Foundation
of the SU-HSE award No. 10-09-0015, and Simons-IUM
fellowship.}
\\[4mm]

{\tt verbit@maths.gla.ac.uk, \ \  verbit@mccme.ru}
\end{center}

{\small 
\hspace{0.15\linewidth}
\begin{minipage}[t]{0.7\linewidth}
{\bf Abstract} \\
Let $(M,I, \omega, \Omega)$ be a nearly K\"ahler
6-manifold, that is, an $SU(3)$-manifold 
with  (3,0)-form $\Omega$ and Hermitian
form $\omega$ which satisfies $d\omega=3\lambda\Re\Omega$,
$d\Im\Omega=-2\lambda\omega^2$, for a non-zero
real constant $\lambda$. We develop an analogue of
the K\"ahler relations on $M$, proving several useful
identities for various intrinsic Laplacians
on $M$. When $M$ is compact, these identities
give powerful results about cohomology of
$M$. We show that harmonic forms on $M$ admit
a Hodge decomposition, and prove that 
$H^{p,q}(M)=0$ unless $p=q$ or $(p=1, q=2)$
or $(p=2, q=1)$.

\end{minipage}
}

{
\small
\tableofcontents
}


\section{Introduction}


\subsection{Nearly K\"ahler 6-manifolds}

Nearly K\"ahler manifolds (also known as 
$K$-spaces or almost Tachibana spaces)
were defined and studied by 
Alfred Gray (\cite{_Gray:1965_}, \cite{_Gray:NK_}, 
\cite{_Gray:weak_holo_}, \cite{_Gray:structure_NK_}) in the general
context of intrinsic torsion of $U(n)$-\-struc\-tures
and weak ho\-lo\-no\-mies. An almost complex Hermitian
manifold $(M,I)$ is called {\bf nearly K\"ahler}
if $\nabla_X(I) X=0$, for any vector field $X$ ($\nabla$ 
denotes the Levi-Civita connection). In other words, the
tensor $\nabla\omega$ must be totally skew-symmetric, for
$\omega$ the Hermitian form on $M$. If
$\nabla_X(\omega)\neq 0$ for any non-zero
vector field $X$, $M$ is called {\bf strictly nearly K\"ahler}.

Using deep results of Kirichenko and Cleyton-Swann
(\cite{_Kirichenko:torsion_NK_}, \cite{_Cleyton_Swann_})
P.-A. Nagy in \cite{_Nagy:splittine_} has shown that 
that any strictly nearly K\"ahler manifold
is locally a product of locally homogeneous manifolds,
strictly nearly K\"ahler 6-manifolds, and twistor
spaces of quaternionic K\"ahler manifolds of
positive Ricci curvature, equipped with
the Eels-Salamon metric (\cite{_Eels_Salamon_}).

These days the term ``nearly K\"ahler'' usually denotes
strictly nearly K\"ahler 6-manifolds. We shall follow
this usage, often omitting ``strictly'' and ``6-dimensional''.
In more recent literature (such as \cite{_MNS:DNK_}),
these objects are called {\bf Gray manifolds}.

For a history of this notion,  a number 
of equivalent definitions and a bibliography
of current work in this field, we refer the reader to 
\cite{_Moroianu_Nagy_Semmelmann:NK_}
and \cite{_Verbitsky:NK_}. 

It is convenient to define nearly K\"ahler 6-manifolds
in terms of differential forms, as follows.

\hfill

\proposition
Let $(M, I, \omega)$ be a Hermitian almost complex 6-manifold.
Then the following conditions are equivalent.
\begin{description}

\item[(i)] The tensor $\nabla_X(I) Y$ is skew-symmetric
with respect to $X$, $Y$, and non-zero.

\item[(ii)] The form $\nabla\omega\in \Lambda^1(M)\otimes \Lambda^2(M)$
is non-zero and totally skew-symmetric.
Notice that in this case, by Cartan's formula,
we have $d\omega = \nabla\omega$.

\item[(iii)] There is $(3,0)$-form $\Omega$ with $|\Omega|=1$, and 
\begin{equation}\label{_NK_defi_forms_Equation_}
\begin{aligned}
d\omega = &  3\lambda \Re\Omega,\\
d \Im \Omega = &  -2 \lambda\omega^2
\end{aligned}
\end{equation}
where $\lambda$ is a non-zero real constant.
\end{description}
{\bf Proof:} \cite{_Gray:NK_}; see also 
\cite{_Baum_etc_}, or
\cite[Theorem 4.2]{_Verbitsky:NK_}.
\endproof

\hfill

\definition
An $SU(3)$-manifold $(M, \omega, \Omega, I)$
is called {\bf nearly K\"ahler} if \eqref{_NK_defi_forms_Equation_}
holds. 

\hfill

The examples of compact nearly K\"ahler manifolds are
scarce; one may hope that the nearly K\"ahler orbifolds
would occur more often.
The results of this paper are stated for the manifolds,
but they are valid for all nearly K\"ahler
orbifolds, with the same proofs.

The most puzzling aspect of nearly K\"ahler geometry
is a complete lack of non-homogeneous 
examples. With the exception of 4 homogeneous 
cases described below (Subsection
\ref{_exa_NK_intro_Subsection_}), no other compact examples of
strictly nearly K\"ahler 6-manifolds are known to exist.

\subsection{Examples of nearly K\"ahler manifolds}
\label{_exa_NK_intro_Subsection_}

Just as the conical singularities of parallel 
$G_2$-manifolds correspond to nearly K\"ahler manifolds,
the conical singularities of $Spin(7)$-\-ma\-ni\-folds
correspond to the so-called ``nearly parallel'' 
$G_2$-manifolds (see \cite{_Ivanov:torsion_G_2_}).
A $G_2$-manifold $(M, \omega)$ is called {\bf nearly parallel}
if $d\omega = c *\omega$, where $c$ is some constant.
The analogy between nearly K\"ahler 6-manifolds and nearly
parallel $G_2$-manifolds is almost perfect.
These manifolds admit a connection with totally
antisymmetric torsion and have weak holonomy
$SU(3)$ and $G_2$ respectively.
N. Hitchin realized nearly K\"ahler 6-manifolds and nearly
parallel $G_2$-manifolds as extrema of a certain 
functional, called {\bf Hitchin functional} by physicists
(see \cite{_Hitchin:stable_}).

However, examples of nearly parallel $G_2$-manifolds
are found in profusion (every 3-Sasakian manifold
is nearly parallel $G_2$), and compact nearly K\"ahler
manifolds are rare.

Only 4 compact examples are known
(see the list below); 
all of them homogeneous.
In \cite{_Butruille_} it was shown that any 
homogeneous nearly K\"ahler 6-manifold 
belongs to this list. 

\begin{enumerate}

\item The 6-dimensional sphere $S^6$.  The 
almost complex structure on $S^6$ is reconstructed
from the octonion action, and the metric is standard.

\item $S^3\times S^3$, with the complex structure
mapping $\xi_i$ to $\xi'_i$,  $\xi_i'$ to $-\xi_i$,
where $\xi_i$, $\xi'_i$, $i=1,2,3$ is a basis
of left invariant 1-forms on the first and the second 
component.

\item Given a self-dual Einstein Riemannian 4-manifold $M$
with positive Einstein constant, one defines
its {\bf twistor space} $\Tw(M)$ as a total space of
a bundle of unit spheres in $\Lambda^2_-(M)$
of anti-self-dual 2-forms. Then $\Tw(M)$ 
has a natural K\"ahler-Einstein structure $(I_+, g)$,
obtained by interpreting unit vectors
in  $\Lambda^2_-(M)$ as  complex structure
operators on $TM$. Changing the sign of 
$I_+$ on $TM$, we obtain an almost complex structure $I_-$
which is also compatible with the metric $g$
(\cite{_Eels_Salamon_}). 
A straightforward computation insures that
$(\Tw(M), I_-, g)$ is nearly K\"ahler
(\cite{_Muskarov_}).

As N.Hitchin proved,
there are only two compact self-dual
Einstein 4-manifolds: $S^4$ and $\C P^2$.
The corresponding twistor spaces are
$\C P^3$ and the flag space $F(1,2)$.
The almost complex structure operator $I_-$ 
induces a nearly K\"ahler structure on these
two symmetric spaces. 

\end{enumerate}

\subsection{Nearly K\"ahler manifolds in geometry and
  physics}

In \cite{_Verbitsky:NK_} it was shown that,
unless a nearly K\"ahler manifold 
$M$ is locally isometric to a 6-sphere,
the almost complex structure on $M$ is uniquely
determined by the metric. In \cite{_Friedrich:S^6_} 
this result was proved for $S^6$ as well. 
Also in \cite{_Verbitsky:NK_} it was shown that
the metric on $M$ is uniquely determined by the almost
complex structure. 

 Denote by $C(M)$ the Riemannian cone
  of $(M, g)$. By definition, the Riemannian cone is a product $\R^{>0}\times
  M$, equipped with a metric $t^2 g + d t^2$,
  where $t$ is a unit parameter of $\R^{>0}$

The definition of nearly K\"ahler manifolds can be
reformulated in terms of Riemannian geometry, as follows.

\hfill

Recall that a
  spinor $\psi$ is called {\bf a Killing spinor} if $\nabla_X\psi =
  \lambda X\cdot \psi$ for all vector fields $X\in TM$ and a
fixed, non-zero real constant $\lambda$. Any manifold
which admits a Killing spinor is Einstein (\cite{_Baum_etc_}).
The following theorem was proven by C. B\"ar.

\hfill

\proposition
Let $(M, g)$ be a Riemannian 6-manifold.
Then $M$ admits a strictly nearly K\"ahler almost complex 
structure if and only if any of the following
equivalent conditions holds.

\begin{description}
\item[(i)] $M$ admits a real Killing spinor.

\item[(ii)] The Riemannian cone $C(M)$ has holonomy
  $G_2$.
\end{description}

{\bf Proof:} \cite{_Bar:Killing_}. 
\endproof

\hfill

For an in-depth study of Killing spinors, with applications
to 6-dimensional geometry, see \cite{_Baum_etc_}.
From (i) it is apparent that a nearly K\"ahler
manifold is Einstein; indeed, only Einstein
manifolds can admit Killing spinors. 

Nearly K\"ahler manifolds appear as the end result of
several important clas\-si\-fi\-cation-\-type problems -
in classification of manifolds admitting a Killing
spinor, in classification of conical singularities
of $G_2$-manifolds, in classification of manifolds
admitting a connection with totally antisymmetric
and parallel torsion (\cite{_Cleyton_Swann_}) 
and so on. These manifolds are even more important
in physics, being solutions of type II B string theory
(\cite{_Friedrich_Ivanov:physics_}).
In that sense, nearly K\"ahler manifolds are just
as important as the usual Calabi-Yau threefolds.

The conical singularities of $G_2$-manifolds and the resulting
nearly K\"ahler geometries also known to have applications
in physics, giving supergravity solutions which are
a product an anti-de Sitter space with an Einstein space 
(see \cite{_Acharya_Jose_Hull_Spence_}). More recently,
the conical singularities of $G_2$-manifolds
arising from nearly K\"ahler geometry were used
to obtain string models with chiral matter fields
(\cite{_Atiyah_Witten_}, \cite{_Acharya_Witten_}).

\subsection{Local structure of nearly K\"ahler 6-manifolds}

Let $(M, I, \omega, \Omega)$ be a nearly K\"ahler manifold.
Since $d\omega = 3\lambda\Re\Omega$, the complex structure
on $M$ is non-integrable; indeed, a differential
of the $(1,1)$-form $\omega$ lies in 
$\Lambda^{3,0}(M)\oplus \Lambda^{0,3}(M)$,
and this is impossible if $(M,I)$ is integrable. 

To fix the notation, we recall some well-known results and
calculations relating the de Rham differential and
the Nijenhuis tensor.

An obstruction to integrability of an 
almost complex structure is given by the 
Nijenhuis tensor, 
\[ N^*:\; T^{1,0}(M)\otimes T^{1,0}(M)\arrow T^{0,1}(M),
\]
mapping a pair of $(1,0)$-vector fields to the
$(0,1)$-part of their commutator. For our purposes,
it is more convenient to deal with its dual,
which we denote by the same letter:
\begin{equation}\label{_Nije_Equation_}
N:\; \Lambda^{0,1}(M)\arrow \Lambda^{2,0}(M).
\end{equation}
{}From Cartan's formula it is apparent that $N$ is
equal to the $(2,-1)$-part of the de Rham differential.
On the other hand, $N$ can be expressed through
$\nabla I$, in a usual way: 
\[
N^*(X,Y) = (\nabla_XI)Y -( \nabla_YI) X
\]
where $X,Y$ are $(1,0)$-vector fields.
On a nearly K\"ahler manifold, $\nabla(I)$ can be
expressed through the 3-form $d\omega=\nabla\omega$.
This gives the following relations (\cite{_Kob_Nomizu_}):
\begin{equation}\label{_Nije_via_d_omega_Equation_}
N^*(X,Y)= d\omega(X, Y, \cdot)^\sharp,
\end{equation}
where $d\omega(X, Y, \cdot)^\sharp$
is a vector field dual to the 1-form $d\omega(X, Y, \cdot)$.
Since $d\omega = 3\lambda\Re\Omega$, the relation 
\eqref{_Nije_via_d_omega_Equation_}
allows one to express $N$ through $\Omega$
and $\omega$. 

Let $\xi_1, \xi_2, \xi_3\in \Lambda^{1,0}(M)$ be an 
orthonormal coframe, such that 
$\Omega= \xi_1\wedge\xi_2\wedge\xi_3$.
Then \eqref{_Nije_via_d_omega_Equation_} gives
\begin{equation}\label{_Nije_in_basis_Equation_}
N(\bar\xi_1) = \lambda\xi_2\wedge \xi_3, \ \ 
N(\bar\xi_2) = -\lambda\xi_1\wedge \xi_3, \ \ 
N(\bar\xi_3) = \lambda\xi_1\wedge \xi_2,
\end{equation}
This calculation is well known; it is
explained in more detail in \cite{_Verbitsky:NK_}.

\subsection{Hodge decomposition of the de Rham 
differential and intrinsic Laplacians}

The results of this paper can be summarized as follows.
Let 
\[
d = d^{2,-1}+ d^{1,0} + d^{0,1}+ d^{-1,2},
\]
be the Hodge decomposition of de Rham differential
(Subsection \ref{_Hodge_on_de_Rham_diffe_Subsection_}).
We use the following notation:
$d^{2,-1}=: N$, $d^{-1, 2}=: \bar N$,
$d^{1,0}=: \6$, $d^{0,1}=: \bar \6$.

The usual K\"ahler identities have a form
``a commutator of some Hodge component of de Rham
differential with the Hodge operator $\Lambda$
is proportional to a Hermitian adjoint of some
other Hodge component of de Rham
differential''. We prove that a similar
set of identities is valid on nearly
K\"ahler manifolds
(\ref{_almost_complex_Kahler_ide_Theorem_},
\ref{_Kodaira_ide_for_N_Proposition_}). These identities
are used to study various intrinsic Laplacians on $M$.
We show that the difference
\begin{equation}\label{_Delta_6-Delta_bar_6_Equation_}
\Delta_\6-\Delta_{\bar\6}= R
\end{equation}
is a scalar operator, acting on $(p,q)$-forms
as $\lambda^2 (p-q) (3-p-q)$ (see
\ref{_diffe_Lapla_6_bar_6_Corollaty_}).
For the de Rham Laplacian $\Delta_d= d d^*+d^*d$, the 
following formula holds:
\begin{equation}
\Delta_d = \Delta_{\6-\bar\6} + \Delta_N + \Delta_{\bar N}
\end{equation}
(see \eqref{_sum_three_Lapla_Equation_}).
This formula is used to study the harmonic
forms on $M$ when $M$ is compact. We show that
$\eta$ is harmonic if and only if all Hodge 
components of $d$ and $d^*$ vanish on $\eta$
(\ref{_Hodge_main_Theorem_}).
This implies that the harmonic forms on
$M$ admit a Hodge decomposition:
\[
{\cal H}^*(M) = \bigoplus {\cal H}^{p,q}(M).
\]
Using \eqref{_Delta_6-Delta_bar_6_Equation_},
we obtain that ${\cal H}^{p,q}(M)=0$
unless $p=q$ or ($q=2, p=1$) or ($q=1, p=2$). We also prove
that all harmonic forms $\eta\in{\cal H}^{p,q}(M)$,
for ($q=2, p=1$) or ($q=1, p=2$) or $p=q=2$ are
{\bf coprimitive}, that is, satisfy
$\eta\wedge \omega=0$, where $\omega$
is the Hermitian  form (\ref{_primitivity_Remark_}).


\section[Algebraic differential operators on the
de Rham algebra]{Algebraic differential operators on the\\
de Rham algebra}
\label{_algebra_diffe_Section_}


The following section is purely algebraic.
We reproduce some elementary results about algebraic
differential operators on graded commutative algebras.
There results are later on used to study the
de Rham superalgebra.

\subsection{Algebraic differential operators: basic properties}

Let $A^*:= \bigoplus^i A_i$ be a graded commutative
ring with unit. For $a\in A_i$, we denote by $L_a$
the operator of multiplication by $a$: $L_a(\eta)=a\eta$.
A supercommutator of two graded endomorphisms
$x, y \in \End(A^*)$ is denoted by
\[
\{x, y\} := xy - (-1)^{\tilde x \tilde y}yx,
\]
where $\tilde x$ denotes the parity of $x$.

Speaking of elements of graded spaces further
on in this section, we shall always mean
{\bf pure} elements, that is, elements of pure
even or pure odd degree. The {\bf parity}
$\tilde x$ is always defined as 1 on odd
elements, and 0 on even elements.

Vectors of pure even degree are called
{\bf even}, and vectors of pure odd degree
are called {\bf odd}. 

A supercommutator of two even endomorphisms, or
an odd and an even endomorphism is equal to their
commutator. We shall sometimes use the usual bracket notation 
$[\cdot, \cdot]$ in this case.

\hfill

\definition\label{_algebra_diffe_o_Definition_}
 The space $D^i(A^*)\subset \End(A^*)$
of {\bf algebraic differential operators of algebraic order $i$} is a graded subspace 
of $\End(A^*)$, which is defined inductively as follows.
\begin{description}

\item[(i)] $D^0(A)$ is a space of $A^*$-linear endomorphisms
of $A^*$, that is, $D^0(A^*)\cong A^*$.

\item[(ii)] $D^{n+1}(A^*)$ is defined as a graded
subspace of $\End(A^*)$ consisting of all 
endomorphisms $\rho\in \End(A^*)$ (even or odd)
which satisfy $\{L_a, \rho\}\in D^{n}(A^*)$, for
all $a\in A$.
\end{description}

This notion was defined by A. Grothendieck.
Using induction, it is easy to check that $D^*(A^*)= \bigcup D^i(A^*)$
is a filtered algebra: 
\begin{equation}\label{_D^*(A)_algebra_Equation_}
D^i(A^*)\cdot D^j(A^*)\subset D^{i+j}(A^*),
\end{equation}
and also
\begin{equation}\label{_D^*(A)_Lie_algebra_Equation_}
\{D^i(A^*), D^j(A^*)\}\subset D^{i+j-1}(A^*).
\end{equation}

\hfill

\definition
Let $\delta:\; A^* \arrow A^*$ be an even or odd endomorphism.
We say that $\delta$ is a {\bf derivation} if
\[
\delta(ab) = \delta(a)b + (-1)^{\tilde a\tilde \delta} a \delta(b),
\]
for any $a, b \in A^*$.

\hfill

Clearly, all derivations of $A^*$ are first order 
algebraic differential operators and vanish on the unit
of $A^*$. The converse is also true: if
$D\subset D^1(A^*)$ is a first order differential
operator, $D(1)=0$, then $D$ is a derivation,
as the following claim implies.

\hfill

\claim\label{_first_order_is_deri_Claim_}
Let $D\in  D^1(A^*)$ be a first order differential
operator. Then \[ D-L_{D(1)}\] is a derivation of $A$.

\hfill

{\bf Proof:} It suffices to prove \ref{_first_order_is_deri_Claim_}
assuming that $D(1)$=0. Let $a, b \in A$ be even or odd elements.
Since $\{ D, L_a\}$ is $A^*$-linear, we have
\begin{align*}
 D(ab) - (-1)^{\tilde a\tilde D} a D(b)= &\{ D, L_a\}(b) 
= \{ D, L_a\}(1) b \\ = &D(a)b + (-1)^{\tilde a\tilde D} a D(1) = D(a)b.
\end{align*}
\endproof

\hfill

\remark\label{_D^1_dete_by_gene_Remark_}
{}From \ref{_first_order_is_deri_Claim_},
it is clear that a first order differential
operator on $A$ is determined by the values
taken on 1 and any set of multiplicative 
generators of $A$.

\hfill

The following claim is also clear. 

\hfill

\claim\label{_D^i_dete_by_gene_Remark_}
Let $D\in \End(A^*)$ be an endomorphism of $A^*$,
and $V$ a set of generators of $A^*$. Assume that
for any $\nu \in V$, we have $\{L_\nu, D\} \in D^i(A^*)$.
Then $D$ is an $(i+1)$-st order algebraic differential
operator on $A$.

\endproof

\subsection{Algebraic differential operators on $\Lambda^*(M)$}

Let $M$ be a smooth manifold, and
$\Lambda^*(M)$ its de Rham algebra. 
It is easily seen that the differential operators
(in the usual sense) and the algebraic differential
operators on $\Lambda^*(M)$ coincide. However,
the ``algebraic order'' of differential operators in the sense
of Grothendieck's definition and in the sense
of the usual definition are different. For instance,
a contraction with a vector field is $C^\infty(M)$-linear, 
hence it has order 0 in the usual sense.
However, the contraction
with a vector field has algebraic order one in
the sense of \ref{_algebra_diffe_o_Definition_}.
Further on, we always use the term ``order'' in
the sense of ``algebraic order'', and not in the
conventional sense.

{}From now till the end of this Appendix,
the manifold $M$ is always assumed to be Riemannian.

\hfill

\claim \label{_contra_1_Claim_}
Let $M$ be a Riemannian manifold,
and $\eta\in \Lambda^1(M)$ a 1-form.
Denote by $\Lambda_\eta$ the metric adjoint to
$L_\eta$, $\Lambda_\eta = -*  L_\eta *$.
Then $\Lambda_\eta$ is a first order differential
operator.

\hfill

{\bf Proof:} 
Clearly, $\Lambda_\eta$ is a contraction with a vector
field $\eta^\sharp$ dual to $\eta$. Then \ref{_contra_1_Claim_}
is clear, because a contraction with a vector field
is clearly a derivation. 
\endproof

\hfill

This claim is a special case of the following proposition,
which is proved independently.

\hfill

\proposition\label{_contra_n_form_Claim_}
Let $M$ be a Riemannian manifold,
and $\eta\in \Lambda^n(M)$ an n-form.
Denote by $\Lambda_\eta$ the metric adjoint to
$L_a$, $\Lambda_\eta = (-1)^{\tilde \eta}  *  L_\eta *$.
Then $\Lambda_\eta$ is a differential
operator of algebraic order $n$.

\hfill

{\bf Proof:} We use the induction on $n$. For $n=0$ everything
is clear. As \ref{_D^i_dete_by_gene_Remark_} implies, to prove
that $\Lambda_\eta\in D^n(\Lambda^*(M))$, we need to show
that 
\begin{equation}\label{_commu_w_L_a_D^n_Equation_}
\{\Lambda_\eta, L_a\} \in D^{n-1}(\Lambda^*(M)),
\end{equation}
for any $a\in \Lambda^0(M), \Lambda^1(M)$.
For $a\in \Lambda^0(M)$, \eqref{_commu_w_L_a_D^n_Equation_}
is clear, because $\Lambda_\eta$ is $C^\infty(M)$-linear,
hence $\{\Lambda_\eta, L_a\}=0$. For $a\in \Lambda^1(M)$,
it is easy to check that
\[
\{\Lambda_\eta, L_a\} = \Lambda_{\eta\cntrct a^\sharp}
\]
where $a^\sharp$ is the dual vector field, and
$\cntrct $ a contraction. The induction statement 
immediately brings \eqref{_commu_w_L_a_D^n_Equation_}.
\endproof

\subsection{An algebraic differential operator and its adjoint}

The main result of this section is the following proposition.

\hfill

\proposition\label{_adjoint_do_D^1_Proposition_}
Let $(M,g)$ be a Riemannian manifold, and 
\[ D:\; \Lambda^*(M)\arrow \Lambda^{*+1}(M)\] a first
order algebraic differential operator. Denote by
$D^*$ its metric adjoint, $D^*=-*D*$. Then
$D^*$ is a second order algebraic differential operator. 

\hfill

{\bf Proof:}
{\bf Step 1:} As follows from \ref{_D^i_dete_by_gene_Remark_},
it suffices to check that
\begin{equation}\label{_seco_commu_line_Equation_}
\{\{D^*, L_a\}, L_b\} \ \ \text{is $\Lambda^*(M)$-linear},
\end{equation}
for all $a, b\in  \Lambda^0(M), \Lambda^1(M)$.

\hfill

{\bf Step 2:}

\hfill

\lemma\label{_first_order_degree_-2_Lemma_}
Let $\Lambda^*(M) \stackrel{D_1}\arrow \Lambda^{*-1}(M)$
be a first order algebraic differential operator
decreasing the degree by 1. Then $D_1=0$.

{\bf Proof:} Follows from \ref{_D^1_dete_by_gene_Remark_}. \endproof

\hfill

{\bf Step 3:} 
 Clearly,
\[ \{\{D^*, L_a\}, L_b\}^* = D_1^*,\]
where $D_1:=\{\{D, \Lambda_a\}, \Lambda_b\}$.
{}From \ref{_contra_1_Claim_} and \eqref{_D^*(A)_Lie_algebra_Equation_},
we find that $D_1$ is an algebraic differential operator
of algebraic order 1
(being a commutator of several first order operators).
When $a, b \in \Lambda^1(M)$, $D_1$ decreases the
degree of a form by 1. By
\ref{_first_order_is_deri_Claim_}, 
$D_1$ is a derivation.
Clearly, a derivation which
vanishes on $\Lambda^0(M)$ is $C^\infty(M)$-linear.
This shows that $D_1$ is $C^\infty(M)$-linear. 

 By \ref{_first_order_degree_-2_Lemma_},
the commutator of $D_1$ with $\Lambda_c$ vanishes, for all
$c\in \Lambda^1(M)$: 
\begin{equation}\label{_commu_D_1_Lambda_Equation_}
\{D_1, \Lambda_c\}=0
\end{equation}

The operator $D_1^*=\{\{D^*, L_a\}, L_b\}$ is $C^{\infty}(M)$-linear
(being adjoint to $D_1$), and commutes with all $L_a$,
as follows from \eqref{_commu_D_1_Lambda_Equation_}.
Therefore, $D_1^*$ is $\Lambda^*(M)$-linear. This proves
\eqref{_seco_commu_line_Equation_} for $a, b \in \Lambda^1(M)$.

\hfill

{\bf Step 4:} Clearly, $L_a=\Lambda_a$ when $a\in \Lambda^0(M)$.
Then 
\begin{equation}\label{_commu_seco_vanishes_Equation_}
\{\{D^*, L_a\}, L_b\}=\{\{D^*, \Lambda_a\}, \Lambda_b\}=
\{\{D, L_a\}, L_b\}^*=0,
\end{equation}
because $D$ is a first order algebraic differential operator.
We proved
 \eqref{_seco_commu_line_Equation_} for $a, b \in \Lambda^0(M)$.

\hfill

{\bf Step 5:} Since the algebra $\Lambda^*(M)$ is graded
commutative, $\{L_a, L_b\}=0$ for all $a, b \in \Lambda^*(M)$.
Using the graded Jacobi identity, we find that
\begin{equation}\label{second_commu_order_Equation}
\{\{D^*, L_a\}, L_b\}= (-1)^{\tilde a \tilde b }\{\{D^*, L_b\}, L_a\},
\end{equation}
for all $a, b$. 
In Steps 3 and 4 we proved 
\eqref{_seco_commu_line_Equation_} for $a, b \in \Lambda^1(M)$,
 $a, b \in \Lambda^0(M)$. By \eqref{second_commu_order_Equation},
to prove \ref{_adjoint_do_D^1_Proposition_}
it remains to show that
$\{\{D^*, L_a\}, L_b\}$ is $\Lambda^*(M)$-linear for
$a\in \Lambda^0(M), b \in \Lambda^1(M)$.

\hfill

{\bf Step 6:} In this case,
\[ \{D^*, L_a\}= 
   \{ D^*, \Lambda_a\} = \{D, L_a\}^*= L_{D(a)}^* = \Lambda_{D(a)}.
\]
Then 
\[ \{\{D^*, L_a\}, L_b\}= \{ \Lambda_{D(a)}, L_b\} = g(D(a), b),\]
because $\Lambda_{D(a)}$ is a contraction with the dual vector
field $D(a)^\sharp$. We proved that $\{\{D^*, L_a\}, L_b\}$
is a scalar function, hence it is $\Lambda^*(M)$-linear.
\ref{_adjoint_do_D^1_Proposition_} is proved. \endproof


\section{K\"ahler identities on nearly K\"ahler manifolds}


\subsection{The operators $\6$, $\bar\6$ on almost complex
  manifolds}
\label{_Hodge_on_de_Rham_diffe_Subsection_}

Let $(M,I)$ be an almost complex manifold, and
$d:\; \Lambda^i(M) \arrow \Lambda^{i+1}(M)$ the de Rham differential.
The Hodge decomposition gives
\begin{equation}\label{_Hodge_on_deRham_Equation_}
d = \bigoplus_{i+j=1} d^{i,j}, \ \ 
 d^{i,j}:\; \Lambda^{p,q}(M)\arrow \Lambda^{p+i,q+j}(M)
\end{equation}
Using the Leibniz identity, we find that the
differential and all its Hodge components
are determined by the values taken
on all vectors in the spaces $\Lambda^0(M)$, $\Lambda^1(M)$,
generating the de Rham algebra. 
On $\Lambda^0(M)$, only $d^{1,0}$, $d^{0,1}$ can be non-zero,
and on $\Lambda^1(M)$ only 
$d^{2,-1}, d^{1,0}, d^{0,1}, d^{-1,2}$ can be non-zero.
Therefore, only 4 components
of \eqref{_Hodge_on_deRham_Equation_} can be possibly non-zero:
\[
d = d^{2,-1}+ d^{1,0} + d^{0,1}+ d^{-1,2},
\]
Since $N:=d^{2,-1}$, $\bar N:=d^{-1,2}$ vanish on
 $\Lambda^0(M)$, these components are $C^\infty(M)$-linear.
In fact, \[ N:\; \Lambda^{0,1}(M)\arrow \Lambda^{2,0}(M)\]
is the Nijenhuis tensor \eqref{_Nije_Equation_} of $(M,I)$, extended to an 
operator on $\Lambda^*(M)$ by the Leibniz rule.
We denote $d^{1,0}$ as $\6:\; \Lambda^{p,q}(M)\arrow \Lambda^{p+1,q}(M)$,
and $d^{0,1}$ as $\bar\6:\; \Lambda^{p,q}(M)\arrow \Lambda^{p,q+1}(M)$.
Decomposing $d^2=0$, we find
\begin{equation}\label{_decompo_d^2_Equation_}
\begin{aligned}
N^2+ \{ N,  \6\}+ & (\{\bar \6, N\} + \6^2) + (\{N, \bar N\} + \{\6, \bar \6\})\\
+&(\{\6, \bar N\} + \bar \6^2)+ \{\bar  N, \bar\6\}  +\bar
N^2 \\
= & d^2=0
\end{aligned}
\end{equation}
where $\{\cdot, \cdot\}$ denotes the supercommutator.
The terms in brackets in \eqref{_decompo_d^2_Equation_}
are  different Hodge components of $d^2$, and
since $d^2=0$, they all vanish:
\begin{equation}\label{_decompo_d^2_vanish_Equation_}
\begin{aligned}
N^2= & \{ N,  \6\}= \{\bar \6, N\} + \6^2 =\\
= & \{N, \bar N\} + \{\6, \bar \6\} = \{\6, \bar N\} + \bar \6^2 \\
 = & \{\bar  N, \bar\6\}  = \bar N^2=0
\end{aligned}
\end{equation}
However, the operators $\6^2$ and
$\bar\6^2$ can be non-zero.

The following
almost complex version of the K\"ahler idenities
is quite useful further on in our study.

\hfill

\theorem \label{_almost_complex_Kahler_ide_Theorem_}
Let $(M,I)$ be an almost complex Hermitian 
manifold, $\omega\in \Lambda^{1,1}(M)$
a Hermitian form, and 
$\Lambda_\omega:\; \Lambda^i(M) \arrow \Lambda^{i-2}(M)$
a Hermitian adjoint to $L_\omega(\eta)= \omega\wedge \eta$.
Consider the operators 
$\6, \bar\6:\; \Lambda^i(M) \arrow \Lambda^{i+1}(M)$
constructed above, and let 
$\6^*, \bar\6^*:\; \Lambda^i(M) \arrow \Lambda^{i-1}(M)$
be the corresponding Hermitian adjoint operators.
Assume that $d\omega\in \Lambda^{3,0}(M)\oplus \Lambda^{0,3}(M)$, 
that is, $\6 \omega=\bar\6\omega=0$.
Then 
\begin{equation}\label{_Kodaira_for_6_bar_6_Equation_}
[\Lambda_\omega, \6] = \1 \bar\6^*, \ \  
[\Lambda_\omega, \bar\6] = -\1\6^*
\end{equation}
and
\begin{equation}\label{_Kodaira_for_6^*_bar_6^*_Equation_}
[L_\omega, \6^*] = \1 \bar\6, \ \ [\Lambda_\omega, \6] = \1 \bar\6^*
\end{equation}
{\bf Proof:} 
To prove \ref{_almost_complex_Kahler_ide_Theorem_},
we use essentially the same argument as used in the
proof of the conventional K\"ahler identities in the situation
when a coordinate approach does not work; see e.g.
 the proof of K\"ahler identities in HKT-geometry, obtained in
\cite{_Verbitsky:HKT_}, and the proof of the K\"ahler
identities in locally conformally hyperk\"ahler geometry,
obtained in \cite{_Verbitsky:LCHK_}. 

The equations \eqref{_Kodaira_for_6_bar_6_Equation_}
and \eqref{_Kodaira_for_6^*_bar_6^*_Equation_}
are Hermitian adjoint, hence equivalent.
The two equations \eqref{_Kodaira_for_6^*_bar_6^*_Equation_}
are complex conjugate, hence they are also equivalent.
To prove \ref{_almost_complex_Kahler_ide_Theorem_}
it is sufficient to prove only one of these equations,
say
\begin{equation}\label{_Kodaira_for_6^*_Equation_}
[L_\omega, \6^*] = \1 \bar\6.
\end{equation}
The proof of such a relation follows a general template,
which is given in  Section \ref{_algebra_diffe_Section_}.
There is an algebraic notion of differential operators
on a graded commutative algebra, due to Grothendieck
(\ref{_algebra_diffe_o_Definition_}).
In this sense, the operators $\6^*$, $\bar \6^*$ are
second order algebraic differential operators on $\Lambda^*(M)$
(see \ref{_adjoint_do_D^1_Proposition_}), and $L_\omega$
is $\Lambda^*(M)$-linear. Then $[L_\omega, \6^*]$
(being a commutator of 0-th and 2-nd order algebraic differential
operators on $\Lambda^*(M)$) is a first order 
algebraic differential operator.  To prove a relation between
first order algebraic differential operators on an algebra, such as 
\eqref{_Kodaira_for_6^*_Equation_},
it suffices to check it on any set of generators 
of this algebra (\ref{_D^1_dete_by_gene_Remark_}). 
To prove \ref{_almost_complex_Kahler_ide_Theorem_}
it remains to show that \eqref{_Kodaira_for_6^*_Equation_}
holds on some set of generators, {\em e.~g.} 
1-forms and 0-forms.

Given a function $f \in C^\infty(M)$, we have
\begin{equation}\label{_kodaira_on_fu_Equation_}
[L_\omega, \6^*] f = -  \6^* (f\omega).
\end{equation}
However, $\6^*=-* \6 *$ and $*(f\omega)= \bar f \omega^{n-1}$,
where $n=\dim_\C M$. Then 
\[ [L_\omega, \6^*] f = * (\6 \bar f \wedge \omega^{n-1})= \1\bar\6 f
\]
because for any $(1,0)$-form
$\eta$ we have $*(\eta \wedge \omega^{n-1})= \1\bar\eta$,
and $\overline{\6\bar f}= \bar\6f$.

It is easy to check that $\6^*$-closed 1-forms
generate the bundle of all 1-forms over $C^\infty(M)$.
Indeed, on 2-forms we have $(\6^*)^2=0$, and therefore,
all $\6^*$-exact 1-forms are $\6^*$-closed.
A local calculation implies that 
$\6^*(\Lambda^2(M))$
generates $\Lambda^1(M)$ over $C^\infty(M)$.

Consider a 1-form $\eta \in \Lambda^1(M)$.
To prove \ref{_almost_complex_Kahler_ide_Theorem_}
it remains to show that
$[L_\omega, \6^*](\eta) = \1 \bar\6\eta$.
Since $\6^*$-closed 1-forms generate
$\Lambda^*(M)$, we may assume that
$\eta$ is $\6^*$-closed. Then
\begin{equation}\label{_commu_L_6_expli_Equation_}
\begin{aligned}
\ [L_\omega, \6^*](\eta) =& -\6^* L_\omega\eta = *\6* (\omega\wedge\eta)\\
= & *\6(\omega^{n-2}\wedge I(\bar\eta)) =   * (\omega^{n-2} \wedge\6 (I\bar\eta)).
\end{aligned}
\end{equation}
Since $\eta$ is $\6^*$-closed, we have $\omega^{n-1} \wedge\6(I\bar\eta)=0$,
and therefore the form $\omega^{n-2} \wedge\6\eta$ is coprimitive
(satisfies $(\omega^{n-2} \wedge\6\eta)\wedge \omega=0$).
Given a coprimitive  $(2n-2)$-form $\alpha= \kappa\wedge\omega^{n-2}$, 
the form $*\alpha$ can be written down explicitly in
terms of $\kappa$: $*\alpha = -I(\bar \kappa)$. 
Then
\[ * (\omega^{n-2} \wedge\6 I(\bar \eta))=- \overline{I \6 I \bar\eta}= \1\bar\6\eta.
\]
Comparing this with \eqref{_commu_L_6_expli_Equation_},
we find that
\[
[L_\omega, \6^*](\eta)=\1\bar\6\eta
\]
We proved \ref{_almost_complex_Kahler_ide_Theorem_}.
\endproof

\subsection{The Nijenhuis operator squared}

Later on, we shall need 
the following useful identity.

\hfill

\proposition\label{_C^2_Proposition_}
Let $(M, I, \omega, \Omega)$ be a nearly K\"ahler
6-manifold, $d\omega=\lambda\Re \Omega$,
and $C:= N+\bar N=d^{2,-1}+d^{-1,2}$ the 
$(2,-1)\oplus(-1,2)$-part of the de Rham differential.
Then the following $C^\infty(M)$-linear maps 
\[ \Lambda^{p,q}(M)\arrow \Lambda^{p+1,q+1}(M)\]
are equal:
\begin{description}
\item[(i)] $C^2$
\item[(ii)] $-\{\6, \bar \6\}$.
\item[(iii)] the scalar operator
  $\1\lambda^2(p-q)L_\omega$,
mapping $\eta \in \Lambda^{p,q}(M)$ to
\[ \1\lambda^2 (p-q) \eta \wedge \omega. \]
\end{description}
{\bf Proof:} The equivalence $C^2=-\{\6, \bar \6\}$
is clear, because the $(1,1)$-part of
$d^2$ is equal to $C^2+ \{\6, \bar \6\}$,
and $d^2=0$ (see \eqref{_decompo_d^2_vanish_Equation_}).
To prove
\begin{equation}\label{_square_of_C_Equation_}
C^2=\1\lambda^2(p-q)L_\omega,
\end{equation}
we notice that both sides of
\eqref{_square_of_C_Equation_}
are differentiations ($C^2$ being a supercommutator
of an odd differentiation with itself). Therefore
it suffices to check \eqref{_square_of_C_Equation_}
only on the generators of $\Lambda^*(M)$, e.g.
on $\Lambda^0(M)$ and $\Lambda^1(M)$.
On $\Lambda^0(M)$, both $C$ and $(p-q)$
vanish, hence \eqref{_square_of_C_Equation_}
is clear. Let us check \eqref{_square_of_C_Equation_}
on $\Lambda^{1,0}(M)$ (a proof of
\eqref{_square_of_C_Equation_} on $\Lambda^{0,1}(M)$
is obtained in the same fashion). Choose an orthonormal
frame $\xi_1, \xi_2, \xi_3\in \Lambda^{1,0}(M)$, in such a way that
\[ \omega = -\1(\xi_1\wedge \bar \xi_1 +\xi_2\wedge \bar
\xi_2 + \xi_3\wedge \bar \xi_3), 
\ \ \ \Omega= \xi_1\wedge \xi_2\wedge \xi_3
\]
Let $\eta$ be a $(0,1)$-form, say, $\eta=\bar\xi_1$
(this assumption is not restrictive, because
both sides of \eqref{_square_of_C_Equation_} are
manifestly $C^\infty(M)$-linear). 
Then $N(\eta)=\lambda\xi_2\wedge \xi_3$, as 
\eqref{_Nije_in_basis_Equation_}
implies.
Similarly, the Leibniz rule and \eqref{_Nije_in_basis_Equation_} give
\begin{equation}\label{_bar_N_N_coordi_Equation_}
\bar N N(\eta)=\lambda^2(\bar \xi_1\wedge \bar \xi_3\wedge 
   \xi_3 + \xi_1\wedge \bar \xi_1\wedge \bar \xi_2)= 
   \1 \lambda^2\eta\wedge\omega.
\end{equation}
On the other hand, 
\begin{equation}\label{_C^2_via_N_equation_}
 C^2(\eta) = (N+\bar N)^2\eta = \bar N N(\eta),
\end{equation}
because $\bar N\eta$ vanishes, being a $(-1,3)$-form,
and $N^2=\bar N^2=0$ as \eqref{_decompo_d^2_Equation_} 
implies. Combining \eqref{_bar_N_N_coordi_Equation_}
and \eqref{_C^2_via_N_equation_}, we obtain 
\eqref{_square_of_C_Equation_}. 
\ref{_C^2_Proposition_} is proved. 
\endproof

\hfill

\corollary \label{_diffe_Lapla_6_bar_6_Corollaty_}
Let $(M,I, \omega, \Omega)$ be a nearly K\"ahler
6-manifold, $d\omega= \lambda \Re\Omega$, and
$\6, \bar\6$ the $(1,0)$- and $(0,1)$-parts
of de Rham differential. Consider the corresponding
Laplacians:
\[
\Delta_\6:= \6\6^*+ \6^*\6, \ \ \ 
\Delta_{\bar\6}:= \bar\6\bar\6^*+ \bar\6^*\bar\6.
\]
Then $\Delta_\6 - \Delta_{\bar\6}= R$,
where $R$ is a scalar operator acting on $(p,q)$-forms
as a multiplication by $\lambda^2 (3-p-q)(p-q)$.

\hfill

{\bf Proof:} As \ref{_C^2_Proposition_} 
 implies, $\{\6, \bar\6\} = -\1 (p-q)\lambda^2 L_\omega$.
It is well known that $H:= [L_\omega, \Lambda_\omega]$
acts on $(p,q)$-forms as a multiplication by $(3-p-q)$
(see e.g. \cite{_Griffi_Harri_}). Then 
\begin{equation}\label{_R_as_commu_Equation_}
\{\Lambda_\omega, \{\6, \bar\6\}\}= \1 R.
\end{equation}
Applying the
graded Jacobi identity and
\ref{_almost_complex_Kahler_ide_Theorem_} 
to \eqref{_R_as_commu_Equation_}, we obtain
\begin{equation}\label{_R_aand_Lapla_Equation_}
\1 R= \{\Lambda_\omega, \{\6, \bar\6\}\}= 
\{ \{\Lambda_\omega, \6\}, \bar\6\} + \{\6,
\{\Lambda_\omega,\bar\6\}\}= \1 \Delta_\6 -\1 \Delta_{\bar\6}.
\end{equation}
This proves \ref{_diffe_Lapla_6_bar_6_Corollaty_}.
\endproof


\section{The de Rham Laplacian via $\Delta_\6$, $\Delta_{\bar
    \6}$, $\Delta_{\6-\bar\6}$}


\subsection{An expression for $\Delta_d$}

Let $M$ be a nearly K\"ahler 6-manifold,
$d = N + \6 + \bar \6 + \bar N$ the Hodge
decomposition of the de Rham differential,
and $\Delta_\6$, $\Delta_{\bar \6}$ the Laplacians
defined above, $\Delta_\6:= \{\6,\6^*\}$, 
$\Delta_{\bar\6}:= \{\bar\6,\bar\6^*\}$.
Denote the usual (Hodge-de Rham) Laplacian by
$\Delta_d=\{d, d^*\}$. 
In addition to \ref{_diffe_Lapla_6_bar_6_Corollaty_},
the following relation between the Laplacians 
can be obtained.

\hfill

\theorem\label{_Laplacians_main_Theorem_}
Let $M$ be a nearly K\"ahler 6-manifold, and
$\Delta_\6$, $\Delta_{\bar\6}$, $\Delta_d$ the
Laplacians considered above. Then 
\begin{equation}\label{_Laplacians_main_theore_Equation_}
\Delta_d = \Delta_\6+ \Delta_{\bar \6} + \Delta_{N+\bar N}
- \{\6, \bar\6^*\} -  \{\bar\6, \6^*\}.
\end{equation}
where $\Delta_{N+\bar N}$ is defined as a supercommutator
of the $C^{\infty}$-linear operator 
$C:=N+\bar N$ and its Hermitian adjoint:
\[
\Delta_{N+\bar N}:= C C^* + C^* C.
\]
The proof of \ref{_Laplacians_main_Theorem_}
takes the rest of this Section.

\subsection{$N= \lambda[L_\Omega, \Lambda_\omega]$}

The following linear-algebraic relation
is used further on in the proof of
\ref{_Laplacians_main_Theorem_}

\hfill

\claim\label{_N_=_[L,Lambda_]_Claim_}
Let $(M,I, \omega, \Omega)$ be a nearly
K\"ahler 6-manifold which satisfies
$d\omega=\lambda\Re\Omega$,
$N$ the $(2, -1)$-part of the
de Rham differential, and $\Lambda_\omega$
the Hodge operator defined above. Then
\begin{equation}\label{_N_=_[L,Lambda_]_Equation_}
\lambda[L_\Omega, \Lambda_\omega]=N,
\end{equation}
where $L_\Omega(\eta):= \Omega\wedge \eta$.

\hfill

{\bf Proof:} As in the proof of 
\ref{_almost_complex_Kahler_ide_Theorem_},
we consider $L_\Omega, \Lambda_\omega$
as algebraic differential operators on the
graded commutative algebra $\Lambda^*(M)$
(see \ref{_algebra_diffe_o_Definition_}). Then $L_\Omega$ is a 0-th order
operator, and $\Lambda_\omega$ a second order operator,
as follows from \ref{_contra_n_form_Claim_}. Therefore, the commutator
$[L_\Omega, \Lambda_\omega]$ is a first
order operator. By \ref{_D^1_dete_by_gene_Remark_} it suffices
to check \eqref{_N_=_[L,Lambda_]_Equation_}
on 0- and 1-forms. This can be done 
by an explicit calculation, using \eqref{_Nije_in_basis_Equation_}.
\endproof

\subsection{Commutator relations for $N$, $\bar N$,
$\6$, $\bar \6$}

\proposition\label{_commu_N_6_Proposition_}
Under the assumptions of \ref{_Laplacians_main_Theorem_},
the following anticommutators vanish:
\begin{equation}\label{_commu_6_and_N_vanish_Equation_}
\{N^*, \bar \6\}=\{\bar N^*, \6\}= \{N, \bar \6^*\}=\{\bar N, \6^*\}=0.
\end{equation}
Moreover, we have
\begin{equation}\label{_commu_6_and_N_via_6_bar_6^*_Equation_}
\begin{aligned}
\{\bar\6^*, \6\} &= -\{ N, \6^*\} = -\{ \bar N^*, \bar\6\}, \\
 \{\6^*, \bar\6\}& = -\{ \bar N,\bar\6^*\}= -\{ N^*, \6^*\},
\end{aligned}
\end{equation}

{\bf Proof:} Clearly, all relations
of \eqref{_commu_6_and_N_vanish_Equation_}
can be obtained by applying the complex
conjugation and taking the Hermitian adjoint of the following relation:
\begin{equation}\label{_commu_6_and_N_one_vanish_Equation_}
\{N^*, \bar\6\}=0.
\end{equation}
Decomposing $d^2=0$ onto Hodge components, we obtain
$\{N, \6\}=0$
(see \eqref{_decompo_d^2_vanish_Equation_}).
By \ref{_N_=_[L,Lambda_]_Claim_}, this is equivalent to 
\begin{equation}\label{_three_commu_L,Lambda,6_Equation_}
\{\{L_\Omega,\Lambda_\omega\}, \6\}=0.
\end{equation}
Clearly, $\6\Omega=0$, hence 
\begin{equation}\label{_L_Omega_6_Equation_}
\{L_\Omega, \6\}=0.
\end{equation}
Applying the graded Jacobi identity to
\eqref{_three_commu_L,Lambda,6_Equation_}
and using \eqref{_L_Omega_6_Equation_},
we obtain 
\begin{equation}\label{_L_Omega_6_triple_Equation_}
0 = \{L_\Omega, \{ \Lambda_{\omega}, \6\}\} = \1\{L_\Omega,\bar\6^*\}.
\end{equation}
Acting on \eqref{_L_Omega_6_triple_Equation_} by 
$\{\Lambda_\omega,\cdot\}$ and using the graded Jacobi
identity, we obtain
\[
0 = \{\Lambda_\omega, \{L_\Omega,\bar\6^*\}\} = -\lambda \{ N,
\bar\6^*\}.
\]
This proves \eqref{_commu_6_and_N_one_vanish_Equation_}
and \eqref{_commu_6_and_N_vanish_Equation_}.

It remains to prove
\eqref{_commu_6_and_N_via_6_bar_6^*_Equation_}.
Taking a Hodge component of $d^2=0$, 
we obtain 
\begin{equation}\label{_6_sqare_Equation_}
\frac 1 2 \{\6,\6\} +\{ N, \bar\6\}=0
\end{equation}
(see \eqref{_decompo_d^2_vanish_Equation_}).
Using the same argument
as gives \eqref{_L_Omega_6_Equation_}, we find
\begin{align*}
\{N, \bar\6\} = & - \lambda^{-1}\{\{\Lambda_\omega,L_\Omega\}, \bar\6\}\\
=&\lambda^{-1}\{\{L_\Omega, \{ \Lambda_\omega, \bar\6\}\}
-\lambda^{-1}\{ \Lambda_\omega, \{L_\Omega, \bar\6\}\}\\
= &-\1\lambda^{-1}\{L_\Omega, \6^*\}.
\end{align*}
Together with \eqref{_6_sqare_Equation_}, this brings
\begin{equation}\label{_6^2_via_L_Omega_Equation_}
\frac 1 2 \{\6,\6\}=\1\lambda^{-1}\{L_\Omega, \6^*\}.
\end{equation}
Acting on \eqref{_6^2_via_L_Omega_Equation_}
with $\{ \Lambda_\omega, \cdot\}$, we obtain
\[
\1\{\6,\bar\6^*\}=\1\lambda^{-1}\{\{ \Lambda_\omega,L_\Omega\},\6^*\}=
-\1\{ N, \6^*\}.
\]
We obtained the first equation of 
\eqref{_commu_6_and_N_via_6_bar_6^*_Equation_}:
\begin{equation}\label{_6,_bar_6^*_final_Equation_}
\{\6,\bar\6^*\}=-\{ N, \6^*\}.
\end{equation}
We proved 
\ref{_commu_N_6_Proposition_}. \endproof

\hfill

\subsection{The Hodge decomposition of the Laplacian}

Now we can finish the proof of 
\ref{_Laplacians_main_Theorem_}. 
Decomposing $d$, $d^*$ onto Hodge components, we obtain
\begin{equation}\label{_Delta_d_decomposed_Equation_}
\begin{aligned}
\Delta_d=& \bigg(\{N^*, \bar \6\} + \{N, \bar \6^*\}+
           \{\bar N^*,  \6\} + \{\bar N,\6^*\}\bigg) \\
+ & \bigg(\{\6^*, \bar \6\} + \{\bar\6^*,  \6\} +\{ N, \6^*\}
           +\{ \bar N, \bar \6^*\}+ \{ N^*, \6\}
           +\{ \bar N^*, \bar \6\}\bigg) \\
+ & \Delta_\6 + \Delta_{\bar\6} + \Delta_{N+\bar N}.
\end{aligned}
\end{equation}
The first term in brackets vanishes by 
\eqref{_commu_6_and_N_vanish_Equation_},
and the second term is equal $-\{\6,\bar\6^*\}-\{\bar\6,\6^*\}$
as \eqref{_commu_6_and_N_via_6_bar_6^*_Equation_} implies. This gives
\begin{equation}\label{_Laplacian_main_Equation_}
\Delta_d = \Delta_\6 + \Delta_{\bar\6} + \Delta_{N+\bar N}
          -\{\6,\bar\6^*\}-\{\bar\6,\6^*\}.
\end{equation}
We proved \ref{_Laplacians_main_Theorem_}. \endproof

\hfill

The relation \eqref{_Laplacian_main_Equation_}
can be rewritten as the following relation between Laplacians.

\hfill

\corollary\label{_Laplacian_via_Delta_6-bar-6_Corollary_}
Under the assumptions of \ref{_Laplacians_main_Theorem_},
denote by $\Delta_{\6-\bar\6}$ the Laplacian
$\{\6-\bar\6, \6^*-\bar\6^*\}$.
Then
\begin{equation}\label{_Laplacian_via_Delta_6-bar-6_Equation_}
\Delta_d =  \Delta_{\6-\bar\6}+ \Delta_{N+\bar N}.
\end{equation}
{\bf Proof:} As \eqref{_Laplacian_main_Equation_}
implies, to prove \eqref{_Laplacian_via_Delta_6-bar-6_Equation_}
we need to show that
\[
\Delta_{\6-\bar\6} = \Delta_\6 + \Delta_{\bar\6} 
          -\{\6,\bar\6^*\}-\{\bar\6,\6^*\}
\]
This is clear.
\endproof


\section{K\"ahler-type identities for $N$, $\bar N$}


Further on, we shall need the following analogue of 
K\"ahler relations
(\ref{_almost_complex_Kahler_ide_Theorem_}),
but for the $C^\infty(M)$-linear ``outer'' parts of the
de Rham differential, $N$ and $\bar N$.

\hfill

\proposition\label{_Kodaira_ide_for_N_Proposition_}
Let $(M, I, \omega, \Omega)$ be a nearly K\"ahler
6-manifold, $N$, $\bar N$ the $(2,-1)$- and $(-1,2)$-parts
of the de Rham differential, $N^*$, $\bar N^*$
their Hermitian adjoints operators, and $\Lambda_\omega$
the Hermitian adjoint to $L_\omega(\eta) :=
\omega\wedge\eta$. Then 
\begin{equation}\label{_commu_N_w_L_Equation_}
\begin{aligned}
\ \ [\Lambda_\omega, N^*] = 2\1 \bar N, \ \ \ &
[\Lambda_\omega, \bar N^*] = -2\1  N.\\
\ \ [L_\omega, N] = 2\1 \bar N^*, \ \ \ &
[L_\omega, \bar N] = -2\1  N^*.
\end{aligned}
\end{equation}
{\bf Proof:} The equalities of \eqref{_commu_N_w_L_Equation_}
are obtained one from another by taking complex
conjugation and Hermitian adjoint, 
hence they are equivalent. Therefore, it suffices to prove
\begin{equation}\label{_commu_N^*_w_L_Equation_}
[L_\omega, N^*] = 2\1 \bar N.
\end{equation}
The proof of this formula follows the same lines as
the proof of \ref{_almost_complex_Kahler_ide_Theorem_}.
Again, both sides of \eqref{_commu_N^*_w_L_Equation_}
are first order algebraic differential operators on the algebra
$\Lambda^*(M)$, in the sense of Grothendieck
(\ref{_algebra_diffe_o_Definition_}). 
Therefore, it suffices to check 
\eqref{_commu_N^*_w_L_Equation_} only on 
0-forms and 1-forms (\ref{_D^1_dete_by_gene_Remark_}). 
On 0-forms, both sides of
\eqref{_commu_N^*_w_L_Equation_} clearly vanish. 
To finish the proof of \eqref{_commu_N^*_w_L_Equation_},
it remains to check
\begin{equation}\label{_commu_N^*_w_L_on_1-forms_Equation}
N^* L_\omega (\eta) = 2\1 \bar N(\eta),
\end{equation}
where $\eta$ is a 1-form. Since both sides of 
\eqref{_commu_N^*_w_L_on_1-forms_Equation} vanish on
$(1,0)$-forms, we may also assume that $\eta \in
\Lambda^{0,1}(M)$.

Let $\xi_1, \xi_2, \xi_2$ be an orthonormal frame
in $\Lambda^{1,0}(M)$, satisfying
\[ \omega = -\1(\xi_1\wedge \bar \xi_1 +\xi_2\wedge \bar
\xi_2 + \xi_3\wedge \bar \xi_3), 
\ \ \ \Omega= \xi_1\wedge \xi_2\wedge \xi_3.
\]
Since both sides of \eqref{_commu_N^*_w_L_Equation_}
are $C^\infty(M)$-linear, we need to prove
\eqref{_commu_N^*_w_L_Equation_} only for
$\eta=\bar\xi_1, \bar\xi_2, \bar\xi_2$.
Assume for example that $\eta=\bar\xi_1$.
Then 
\[ L_\omega\eta = - \1 \bar \xi_1 \wedge 
(\xi_2 \wedge \bar\xi_2 + \xi_3\wedge \bar\xi_3).
\]
Using \eqref{_Nije_in_basis_Equation_}, we obtain
\[
N^* L_\omega\eta =-\1 \xi_3\wedge  \xi_2 
+\1 \xi_2\wedge  \xi_3 = 2 \1 \xi_2\wedge\xi_3= 
  2\1\: \bar N(\eta).
\]
This proves \eqref{_commu_N^*_w_L_Equation_}.
\ref{_Kodaira_ide_for_N_Proposition_} is proved.
\endproof

\hfill

 \ref{_Kodaira_ide_for_N_Proposition_} 
is used in this paper only once,  to obtain the following 
corollary.

\hfill

\corollary\label{_N_barN^*_commute_Corollary_}
Let $M$ be a nearly K\"ahler 6-manifold, 
$N$, $\bar N$ the $(2,-1)$- and $(-1,2)$-parts
of the de Rham differential, $N^*$, $\bar N^*$
their Hermitian adjoint operators, and
$\Delta_N$, $\Delta_{\bar N}$, $\Delta_{N+\bar N}$
the corresponding Laplacians. Then 
\[
\Delta_{N+\bar N}=\Delta_N+\Delta_{\bar N}
\]
{\bf Proof:}
Clearly, we have
\[
\Delta_{N+\bar N}=\Delta_N+\Delta_{\bar N}+ \{N, \bar N^*\} + \{\bar N,  N^*\}.
\]
Therefore, to prove \ref{_N_barN^*_commute_Corollary_},
it suffices to show that
\[ 
 \{N, \bar N^*\}=0,  \ \  \{\bar N,  N^*\}=0.
\]
One of these equations is obtained from another
by complex conjugation; therefore, they are equivalent.
Let us prove, for instance, $\{N, \bar N^*\}=0$. As 
follows from \ref{_Kodaira_ide_for_N_Proposition_},
\begin{equation}\label{_N,bar_N^*_commu_Equation_}
\{N, \bar N^*\} = - \frac \1 2 \{ N, \{ \Lambda_\omega, N\}\}.
\end{equation}
However, $\{N, N\}=0$ as follows from 
\eqref{_decompo_d^2_vanish_Equation_}.
Using the graded Jacobi identity, we obtain
\[ 0 = \{ \Lambda_\omega,\{  N, N\} \} =
       \{ \{\Lambda_\omega,  N\}, N \}+ \{ N, \{\Lambda_\omega,  N\}\}
   = 2 \{ N, \{\Lambda_\omega,  N\}\}.
\]
Therefore, \eqref{_N,bar_N^*_commu_Equation_} implies 
$\{N, \bar N^*\}=0$. This proves \ref{_N_barN^*_commute_Corollary_}.
\endproof

\hfill

\remark
{}From \ref{_N_barN^*_commute_Corollary_}, 
\ref{_Laplacian_via_Delta_6-bar-6_Corollary_}
 and \ref{_Laplacians_main_Theorem_}, we infer
that 
\begin{equation}
\Delta_d= \Delta_\6+\Delta_{\bar\6}+ \Delta_N + \Delta_{\bar N}-
\{\6^*,\bar\6\} - \{\6,\bar\6^*\}
\end{equation}
and 
\begin{equation}\label{_sum_three_Lapla_Equation_}
\Delta_d= \Delta_{\6-\bar\6}+ \Delta_N + \Delta_{\bar N}
\end{equation}


\section{Harmonic forms on nearly K\"ahler manifolds}


\subsection{Harmonic forms and the $\6$, $\bar\6$-Laplacians}

For harmonic forms on a compact nearly K\"ahler manifold,
the relation \eqref{_Laplacians_main_theore_Equation_} of
\ref{_Laplacians_main_Theorem_} can be
strengthened significantly.

\hfill

\theorem\label{_harmo_on_compa_NK_main_Theorem_}
Let $M$ be a compact nearly K\"ahler 6-manifold,
and $\eta$ a differential form. Then $\eta$ is harmonic if and
only if 
\begin{equation}\label{_eta_harmo_for_all_Lapla_Equation_}
\eta\in \ker \Delta_\6 \cap \ker \Delta_{\bar\6}\cap \ker \Delta_{N+\bar N}.
\end{equation}
{\bf Proof:} The ``if'' part is clear; indeed, if
\eqref{_eta_harmo_for_all_Lapla_Equation_} is true,
then $\6\eta =\bar\6\eta=\6^*\eta=\bar\6^*\eta=0$,
and $\Delta_d\eta=0$ by \ref{_Laplacians_main_Theorem_}.

As \ref{_Laplacian_via_Delta_6-bar-6_Corollary_} implies,
\begin{equation}\label{_Lapla_implies_2_others_Equation_}
\Delta_d\eta=0 \Leftrightarrow 
\bigg( \Delta_{\6-\bar\6}\eta=0, \ \text{and} \ \ 
\Delta_{N+\bar N}\eta=0\bigg).
\end{equation}  
Therefore, for any $\Delta_d$-harmonic
form $\eta$, we have $\Delta_{\6-\bar\6}\eta=0$, that is,
$(\6-\bar\6)\eta=0$ and $(\6^*-\bar\6^*)\eta=0$.
Moreover, since $d\eta=(N+\bar N)\eta=0$,
$\6+\bar\6= d-N-\bar N$ vanishes on $\eta$ as well.
Substracting from $(\6+\bar\6)\eta=0$
the relation $(\6-\bar\6)\eta=0$, we obtain that
$\bar\6\eta=0$.  In a similar way one proves the whole set
of equations
\[ \6\eta =\bar\6\eta=\6^*\eta=\bar\6^*\eta=0
\]
This gives \eqref{_eta_harmo_for_all_Lapla_Equation_}.
\ref{_harmo_on_compa_NK_main_Theorem_} is proved.
\endproof

\subsection{The Hodge decomposition on cohomology}

The main result of this paper is an immediate corollary of
\ref{_harmo_on_compa_NK_main_Theorem_}.

\hfill

\theorem\label{_Hodge_main_Theorem_}
Let $M$ be a compact nearly K\"ahler 6-manifold,
and ${\cal H}^i(M)$ the space of harmonic $i$-forms on $M$.
Then ${\cal H}^i(M)$ is a direct sum of harmonic forms
of pure Hodge type:
\begin{equation}\label{_Hodge_on_diffe_forms_Equation_}
{\cal H}^i(M)=\bigoplus\limits_{i=p+q}{\cal H}^{p,q}(M).
\end{equation}
Moreover, ${\cal H}^{p,q}(M)=0$ unless $p=q$
or $(p=2, q=1)$ or $(q=1, p=2)$.

\hfill

{\bf Proof:} 
{}From \ref{_harmo_on_compa_NK_main_Theorem_}, we find that a form $\eta$ 
is harmonic if and only if $\6\eta=\bar\6\eta=\6^*\eta=\bar\6^*\eta=0$
and $\Delta_{N+\bar N}\eta=0$. From \ref{_N_barN^*_commute_Corollary_},
we find that the latter equation is equivalent to 
$N\eta= \bar N\eta=N^*\eta = \bar N^*\eta=0$.
We find that a form $\eta$ is harmonic if and only if 
all the Hodge components of $d$, $d^*$ vanish
on $\eta$:
\begin{equation}\label{_harmonic_all_compo_vanish_Equation_}
\6\eta=\bar\6\eta=\6^*\eta=\bar\6^*\eta=
N\eta= \bar N\eta=N^*\eta = \bar N^*\eta=0.
\end{equation}
Therefore, all Hodge components of
$\eta$ also satisfy \eqref{_harmonic_all_compo_vanish_Equation_}.
This implies that these components are also harmonic. 
We proved \eqref{_Hodge_on_diffe_forms_Equation_}.

To prove that ${\cal H}^{p,q}(M)$ vanishes unless $p=q$
or $(p=2, q=1)$ or $(q=1, p=2)$, we use \ref{_diffe_Lapla_6_bar_6_Corollaty_}. 
Let $\eta$ be a non-zero harmonic $(p,q)$-form. Then 
the scalar operator $R=\Delta_\6-\Delta_{\bar\6}$
vanishes on $\eta$, $R= \lambda^2(3-p-q)(p-q)$.
Therefore, either $p=q$ or $p+q=3$. We obtain immediately
that $p=q$ or $(p=2, q=1)$ or $(q=1, p=2)$ or
$(p=3, q=0)$ or $(q=0, p=3)$. The last two cases are
impossible: on $\Lambda^{3,0}(M)$, $\Lambda^{0,3}(M)$,
the operator $N+\bar N$ is clearly injective
(see \eqref{_Nije_in_basis_Equation_}), 
hence it cannot vanish; however,
by \eqref{_harmonic_all_compo_vanish_Equation_}
we have $N+\bar N(\eta)=0$. This proves 
\ref{_Hodge_main_Theorem_}. \endproof

\hfill

\remark
The middle cohomology of a compact nearly K\"ahler 6-manifold
is remarkably similar to the middle cohomology of a K\"ahler manifold.
In particular, the intermediate Jacobian $T:=H^{2,1}(M)/H^3(M, \Z)$
is well defined in this case as well. As in the K\"ahler case,
$T$ is a compact complex torus, and we have a
pseudoholomorphic map
$S\arrow T$ from the space of pseudoholomorphic rational
curves on $M$ to the intermediate Jacobian.

\hfill

\remark\label{_primitivity_Remark_}
All harmonic forms $\eta\in {\cal H}^{p,q}(M)$,
$(p=2, q=1)$ or $(q=1, p=2)$ are {\bf primitive}
and {\bf coprimitive},
that is, satisfy $L_\omega(\eta)=\Lambda_\omega(\eta)=0$.
Indeed, $\{ N, \bar N\}$ vanishes on $\eta$
as follows from
\eqref{_harmonic_all_compo_vanish_Equation_}.
However, \ref{_C^2_Proposition_} implies that
$\{ N, \bar N\}(\eta) =(p-q)\omega\wedge \eta$,
hence $\omega\wedge \eta=0$. 

Similarly, all harmonic forms $\eta\in {\cal H}^{p,p}(M)$
are primitive for $p=1$ and comprimitive for $p=2$.
This  is implied directly by
\eqref{_harmonic_all_compo_vanish_Equation_}
and the local expression for $N$ given in
\eqref{_Nije_in_basis_Equation_}. For instance,
for a $(1,1)$-form $\eta$, we have
$N(\eta) = \Lambda_\omega (\eta)\cdot \Omega$,
hence $\eta$ is primitive if $N(\eta)=0$.


\section{Appendix: Hodge theory on orbifolds}


In the appendinx, we explain how the results of this paper
can be applied to compact nearly K\"ahler orbifolds.

Recall that {\bf an orbifold} is a topological
space equipped with an atlas of local charts, which
are isomorphic to $\R^n /G$, where $G$ is a finite
group acting faithfully and smoothly, and with
all the gluing maps smooth and compatible with the group action.

The differential forms on an orbifold are defined
in local charts as $G$-invariant differential forms
on $\R^n$, which are compatible with the gluing maps.

The de Rham algebra and de Rham cohomology are defined
literally in the same way as for manifolds, and
they are equal to the singular cohomology. 
This was first observed by I. Satake, who defined
the orbifolds in 1950-ies and called them
``V-manifolds'' (\cite{_Satake:V_manifold_},
\cite{_Satake:Gauss-Bonnet_}). Since Satake, 
all the usual constructions of smooth topology,
such as the Chern-Weil theory of characteristic
classes, Atiyah-Singer index formula, signature
theorem and Riemann-Roch-Grothendieck theorem,
were generalized for the orbifolds
(see e.g. \cite{_Kawasaki:signature_}, 
\cite{_Kawasaki:RR_}, \cite{_Kawasaki:index_}).

The Hodge theory identifies harmonic forms with the
de Rham cohomology, using the closedness of the image
of the de Rham differential on a compact manifold.
The basic machinery here works for the orbifolds as
well as in the smooth case (\cite{_Kawasaki:index_}).

The results of the present paper, such as 
\ref{_Hodge_main_Theorem_}, are 
obtained by application of local formulas for
Laplacians to the global statements about cohomology,
using the identification of cohomology and
harmonic forms. These arguments are transferred
to the orbifold case word by word.

The vanishing theorems, such as
\ref{_harmo_on_compa_NK_main_Theorem_}, 
are obtained by showing that
a difference of certain second order operators is positive, which implies
that a kernel of one of these operators lies inside a kernel
of another. Here
the vanishing arguments are also translated to the
orbifold case, without any difficulty.

\hfill

{\bf Acknowledgements:} I am grateful to Robert Bryant
and Paul-Andi Nagy for enlightening correspondence.
P.-A. Nagy also suggested to add \ref{_primitivity_Remark_}.

\hfill

{\scriptsize

\hfill

\small

\noindent {\sc Misha Verbitsky\\
Laboratory of Algebraic Geometry,
Faculty of Mathematics, NRU HSE,
7 Vavilova Str. Moscow, Russia }\\
\tt verbit@maths.gla.ac.uk, \ \  verbit@mccme.ru}

\end{document}